\documentclass[12pt]{amsart}

\usepackage{amssymb,latexsym}
\usepackage{amsmath,amsfonts,amsthm,amscd,amsxtra,enumerate}
\usepackage{mathtools}
\usepackage[all]{xy}
\usepackage{amsmath,latexsym,amssymb, enumerate,  amsthm, epsfig, hyperref}

\usepackage[margin=1in]{geometry}
\usepackage{epsfig}

\newtheorem{theorem}{Theorem}[section]
\newtheorem{lemma}[theorem]{Lemma}
\newtheorem{corollary}[theorem]{Corollary}
\newtheorem{definition}[theorem]{Definition}
\newtheorem{example}[theorem]{Example}
\newtheorem{remark}[theorem]{Remark}
\newtheorem{proposition}[theorem]{Proposition}

\newtheorem{question}[theorem]{Question}

\newtheorem{setup}[theorem]{Setup}

\newcommand{\p}{\mathbf{P}}
\newcommand{\sk}{{\ensuremath{\sf k }}}
\newcommand{\m}{\ensuremath{\mathfrak m}}
\newcommand{\ZZ}{\ensuremath{\mathbb{Z}}}
\newcommand{\QQ}{\ensuremath{\mathbb{Q}}}

\newcommand{\Y}{\ensuremath{\mathbf{Y}}}
\newcommand{\Z}{\ensuremath{\mathbf{Z}}}
\newcommand{\y}{\ensuremath{\mathbf{y}}}
\newcommand{\z}{\ensuremath{\mathbf{z}}}

\newcommand{\eR}{\ensuremath{\varepsilon_R}}
\newcommand{\eS}{\ensuremath{\varepsilon_S}}
\newcommand{\lm}{\ensuremath{\lambda}}
\newcommand{\ov}{\ensuremath{\overline}}

\renewcommand{\ll}{{\ensuremath{\ell\ell}}}

\DeclareMathOperator{\ann}{ann}
\DeclareMathOperator{\Tor}{Tor}
\DeclareMathOperator{\soc}{soc}

\DeclareMathOperator{\edim}{edim}
\DeclareMathOperator{\type}{type}
\DeclareMathOperator{\gr}{gr}

\begin{document}

\title{\textbf{Decomposing Gorenstein Rings as Connected Sums}}

\author[H. Ananthnarayan]{H. Ananthnarayan}
\address{Department of Mathematics, I.I.T. Bombay, Powai, Mumbai 400076.}
\email{ananth@math.iitb.ac.in}

\author[Ela Celikbas]{Ela Celikbas}
\address{Department of Mathematics, West Virginia University, Morgantown, WV 26506.}
\email{ela.celikbas@math.wvu.edu}

\author[Jai Laxmi]{Jai Laxmi}
\address{Department of Mathematics, I.I.T. Bombay, Powai, Mumbai 400076.}
\email{jailaxmi@math.iitb.ac.in}

\author[Z. Yang]{Zheng Yang}
\address{Department of Mathematics, Miami University, Oxford, OH 45056.}
\email{yangz15@miamioh.edu}



\thanks{Z. Yang was partially supported by NSF grant DMS-1103176.}

\begin{abstract}
In 2012, Ananthnarayan, Avramov and Moore
give a new construction of Gorenstein rings from two Gorenstein local rings, called their connected sum.
Given a Gorenstein ring, one would like to know whether it decomposes as a connected sum and if so,
what are its components. We answer these questions in the Artinian case 
and investigate conditions on the ring which force it to be indecomposable as a connected sum. We further
give a characterization for Gorenstein Artin local rings to be decomposable as connected sums, and as a consequence, obtain results about its
Poincar$\acute{\text{e}}$ series and minimal number of generators of its
defining ideal. Finally, we show that the indecomposable components appearing in the connected sum decomposition are unique up to isomorphism.
\end{abstract}

\maketitle

\section*{Introduction}

The main object of study in this paper is a construction of Gorenstein rings, called a connected sum, defined by 
Ananthnarayan, Avramov and Moore in \cite{AAM}. Given Cohen-Macaulay local rings $R$, $S$ and $\sk$
of the same dimension, and ring homomorphisms 
$R \overset{\varepsilon_R}\longrightarrow \sk \overset{\varepsilon_S}\longleftarrow S$,
the authors consider the fibre product (or pullback) $R\times_\sk S = \{(r,s) \in R \times S: \eR(r) = \eS(s)\}$ 
and define a connected sum of $R$ and $S$ over $\sk$ as an appropriate quotient of $R\times_\sk S$. They prove that 
when $R$ and $S$ are Gorenstein, a connected sum is also a Gorenstein local ring of the same dimension. 

In this paper, we focus on connected sums over a field in the Artinian case, i.e., when $R$ and $S$ are 
Gorenstein Artinian local rings and $\sk$ is their common residue field. These objects have been studied 
from different perspectives by various authors starting with Sah (cf. \cite{Sah}) in the graded case and, 
in the local case, by Lescot (see Remark \ref{CS}(d)). A topologically influenced version was also studied 
by Smith and Stong (cf. \cite[Section 4]{SS}), and quite a few authors approach this area via 
{\it Macaulay's inverse systems}, (e.g., see \cite{BBKT}). Completely different techniques are used in this 
article: we look at intrinsic properties of the ring and its defining ideal in Section \ref{Indec}.  

A natural question is: Given a Gorenstein Artin local ring $Q$, can it be decomposed as a connected sum?
It is known (see \cite[8.3]{AAM}) that if $Q$ has embedding dimension at least 3, and decomposes non-trivially
as a connected sum over $\sk$, then
$Q$ cannot be a {\it complete intersection}. In the equicharacteristic case, the question of decomposability has been studied from a geometric point of view by Smith and Stong (cf. \cite[Section 4]{SS}) for {\it projective bundle ideals}. This question has also been studied via inverse systems using polynomials that are {\it direct sums} and corresponding {\it apolar} Gorenstein algebras, by Buczy\'nska et al in \cite{BBKT}  (see Remark \ref{apolar}) .

In Section \ref{Indec}, we give conditions for the indecomposibility as connected sum of a Gorenstein Artin local ring $Q$ in terms of its Hilbert function (Theorem \ref{HilbertCoeff2}), and in terms of the minimal number of generators of its defining ideal (Theorem \ref{mingensTheorem}). In particular,
one can see that if $Q$ is {\it compressed} with Loewy length at least 4, then it is
indecomposable as a connected sum (see Corollary \ref{compressed}).

We then identify necessary conditions for $Q$ to be a connected sum
in terms of its defining ideal in Proposition \ref{MainProp}. We use this proposition and Remark \ref{CS}(\ref{CS1}) to give equivalent conditions for a
Gorenstein Artin local ring to be a connected sum over its residue field, which is the content of Theorem \ref{MainThm}. 

A second question is: If $Q$ is a connected sum, what are its components? Remark \ref{CS}(\ref{CS1}), together with Proposition \ref{Prop2}, relates the defining ideals of a connected sum with the defining ideal of its components. Furthermore, in the equicharacteristic case, this allows
us to find the components computationally (see Remark \ref{Prop3}). We then give two applications of the main theorem. If $Q$ is a given Gorenstein Artin $\sk$-algebra, we give a condition on $Q/\soc(Q)$ which forces $Q$ to be decomposable as a connected sum. In the second application, if $Q$ is decomposable as a connected sum, we identify some conditions under which a Gorenstein quotient of $Q$ is also decomposable.

One can also ask whether such a decomposition into connected sums over $\sk$ is unique. We answer this question in the graded case, by proving that the indecomposable components appearing in the connected sum decomposition are unique up to isomorphism in Section \ref{Unique} (see Theorem \ref{3.5}). 

The first two sections contain results regarding the main tools used in the rest of the paper.
In Section \ref{Preliminaries}, we collect some properties of associated graded rings, Cohen presentations and Poincar$\acute{\text{e}}$ series.
Section \ref{FC} contains some basic observations about fibre products and connected sums, including their interactions
with the objects introduced in Section 1. We then develop some technical results regarding Cohen presentations, in particular, we give the Cohen presentation of a fibre product of complete regular local rings. In particular, Proposition \ref{CP2} gives the relation between defining ideals of a fibre product ring and its components. These results are used in the rest of the article, to get the Cohen presentations of fibre products and connected sums. The results in Section 1 are well-known and Section 2 includes
some of our basic observations, and other known results rephrased in our notation. 

The computer algebra package Macaulay2 has been used for computations in some of the examples in this article.

\section{Preliminaries}\label{Preliminaries}
\subsection{Notation}\hfill{}
\begin{enumerate}[{\rm a)}]
\item  For positive integers $m$ and $n$, $\Y$ and $\Z$ denote the sets $\{Y_1, \ldots, Y_m\}$ and
$\{Z_1,\ldots,Z_n\}$ respectively, and $\Y \cdot \Z$ denotes $\{Y_iZ_j: 1 \leq i \leq m, 1 \leq j \leq n\}$.

\item If $\sk$ is a field, a {\it graded $\sk$-algebra} $G$ is a graded ring $G = \oplus_{i \geq 0} G_i$ with
$G_0 = \sk$. It has a unique homogeneous maximal ideal, $G_+ = \oplus_{i \geq 1} G_i$. We say $G$ is {\it standard graded} if $G_+$ is generated by $G_1$.

\item  If $T$ is a local ring, $\m_T$ denotes the maximal ideal of $T$. Furthermore, for a $T$-module $M$, $\lambda(M)$ and $\mu(M)$ respectively denote the {\it length} and the
{\it minimal number of generators} of $M$ as a $T$-module.

\item Let $(T,\m,\sk)$ be an Artinian local ring. Then $\edim(T)$
denotes the {\it embedding dimension} of $T$ which is equal to $\mu(\m)$.
The {\it socle} of $T$ is $\soc(T) = \ann_T(\m)$. Moreover, the {\it type} of $T$ is $\type(T) = \dim_{\sk}(\soc(T))$,
and the \emph{Loewy length} of $T$ is $\ll(T) = \max\{n: \m^n \neq 0\}$.\footnote{If $T$ is also Gorenstein,
its Loewy length is also referred to as {\it socle degree} in the literature.}

Observe that $T$ is not a field if and only if $\ll(T) \geq 1$. Furthermore, if $T$ is Gorenstein, then $\soc(T) \subset \m^2$ if and only if $\ll(T) \geq 2$.
\end{enumerate}

\subsection{Associated graded rings}
\begin{definition}{\rm Let $(T,\m,\sk)$ be a Noetherian local ring.\begin{enumerate}[{\rm a)}]
\item The {\it graded ring associated to the maximal ideal} $\m$ of
$P$, denoted $\gr_{\m}(T)$, {\rm (}or simply $\gr(T)${\rm)}, is defined as $\gr(T) \simeq \oplus_{i=0}^{\infty} (\m^i/\m^{i+1})$.

\item Let $G = \gr(T)$. We define the {\it Hilbert function} of $T$ as $H_T(i) = \dim_{\sk}(G_i)$ for $i \geq 0$.

\item When $T$ is Artinian, we write the Hilbert function of $T$
as $H_T=(H_T(0),\ldots,H_T(s))$, where $\ll(T) = s$.

Furthermore, if $\ T$ is Gorenstein, we say that $T$ is \emph{short} if $\m^4_T = 0$, and $T$ is \emph{stretched} if $\ \m_T^2$ is principal, i.e.,
$H_T(i) = 1$ for $2 \leq i \leq s$.
\end{enumerate}
}\end{definition}

\begin{remark}\label{GR}{\rm With notation as above, let $G = \gr(T)$ and
$G_{\geq n} = \bigoplus_{i=n}^{\infty} \m^i/\m^{i+1}$ for $n \geq 0$.

\begin{enumerate}[a)]
\item For each $n \geq 0$, $G_{\geq n}$ is the $n$th power of the homogeneous maximal ideal $G_+$ of $G$ and a minimal
generating set of $G_{\geq n}$ lifts to a minimal generating set of $\m^n$.

In particular, if $T$ is Artinian local, then so is $\gr(T)$. Furthermore, $\lambda(T)  = \lambda(\gr(T))$ and $\ll(T) = \ll(\gr(T))$.

\item For each $x \in T \setminus \{0\}$, there exists a unique non-negative integer $i$ such that $x \in \m^i \setminus \m^{i+1}$. The
{\it initial form} of $x$ is the element $x^* \in G$ of degree $i$ that is the image of $x$ in $\m^i/\m^{i+1}$.

\item For an ideal $K $ of $T$, $K^*$ denotes the ideal of $G$ defined by $\langle x^*: x \in K \rangle$.
Note that, if $R \simeq T/K$, then $\gr(R) \simeq G/K^*$.
\end{enumerate}
}\end{remark}

\subsection{Cohen Presentations and Poincar$\mathbf{\acute{\text{e}}}$ Series}

\begin{definition}
{\rm Let $T$ be a local ring. We say that $\widetilde T / I_T$ is a {\it Cohen presentation} of $T$ if $(\widetilde T, \m_{\widetilde T}, \sk)$ is a complete regular local ring and $I_T \subset \m_{\widetilde T}^2$ is an ideal in $\widetilde T$ such that $T \simeq \widetilde T / I_T$.
}\end{definition}

\begin{remark}\hfill{}{\rm

\begin{enumerate}[a)]
\item By Cohen's Structure Theorem, every complete Noetherian local ring has a Cohen presentation. 
\item Let $(\widetilde T, \m_{\widetilde T}, \sk)$ be a complete regular local ring, $I$ be an ideal in $\widetilde T$. Set $T =  \widetilde T / I$. Then $\widetilde T / I$ is a Cohen presentation of $T$ if and only if $\edim(\widetilde T) = \edim(T)$.
\end{enumerate}
}\end{remark}

\begin{definition}\label{PSdef}{\rm For a local ring $(T, \m, \sk)$, the {\it Poincar$\acute{\text{e}}$ series} of $T$, is the formal power series
$$\p^T(t)=\sum_{i\geq 0} \beta_i^Tt^i, \;\;\; \text{with}\;\; \beta_i^T=\dim_{\sk}\left( \Tor_i^T(\sk, \sk)\right).$$
}\end{definition}

\begin{remark}[Minimal Number of Generators]\label{PS2}{\rm
Let $(T,\m,\sk)$ be a complete Noetherian local ring with $\edim(T) = d$, and $T = \widetilde T/I_T$ be a Cohen presentation.
By \cite[7.1.5]{Av}, we have
$$\beta_1^T =d\quad\text{ and }
\quad\mu(I_T) =\beta_2^T - \binom{\beta_1^T}{2} = \beta_2^T - \binom{d}{2}.$$

}\end{remark}

Next we list some properties of the Poincar$\acute{\text{e}}$ series of a Gorenstein Artin local ring.

\begin{remark}[Poincar$\mathbf{\acute{\text{e}}}$ Series of Gorenstein Rings]\label{PS} \hfill{}\\{\rm
Let $(T,\m,\sk)$ be a Gorenstein Artin local ring and $\overline T$ represent the quotient $T/\soc(T)$.

\begin{enumerate}[a)]
\item If $\edim(T) \geq 2$, then
$[\p^{T}(t)]^{-1} = [\p^{\overline T}(t)]^{-1} + t^2$,  by \cite[Thm. 2]{LA}.

\item Let $\ll(T) = 2$. Then $\p^{T}(t) = (1- t)^{-1}$ if $\edim(T) = 1$, and if $\edim(T) \geq 2$, then $\p^{T}(t) = (1 - nt + t^2)^{-1}$
 (for example, by (a), since $\ll(\overline{T}) = 1$).
\end{enumerate}
}\end{remark}

\section{Fibre Products and Connected Sums}\label{FC}

In this section, we see the definition and some basic properties of fibre products and connected sums (Remarks \ref{FP} and \ref{CS}). Some of the observations in these remarks can be found in \cite{AAM} and \cite[Chapter 4]{HAthesis}, we present them here for the sake of completeness.

\subsection{Fibre Products}
\begin{definition}\label{FPdef}{\rm
Let $(R,\m_R,\sk)$ and $(S,\m_S,\sk)$ be local rings. The {\it fibre product} $R$ and $S$ over $\sk$ is the ring
$R\times_{\sk} S = \{(r,s) \in R \times S: \pi_R(r) = \pi_S(s)\}$,
where $\pi_R$ and $\pi_S$ are the natural projections from $R$ and $S$ respectively onto $\sk$.
}\end{definition}

\begin{remark}\label{FP}{\rm
With the notation as in Definition \ref{FPdef}, set $P = R\times_{\sk} S$. 
\begin{enumerate}[a)] 
\item (Trivial fibre product). Every ring is trivially a fibre product over its residue field. Indeed, if $S \simeq \sk$, then 
$P = R \times_{\sk} \sk \simeq R$. 

\item We have a short exact sequence of $P$-modules:\\ $0 \rightarrow P \overset{\epsilon}\rightarrow R \oplus S \overset{\pi}\rightarrow \sk \rightarrow 0$, where $\epsilon(r,s) = (r,s)$ and $\pi(r,s) = \pi_R(r) - \pi_S(s)$. 

In particular, $\widehat P \simeq \widehat R \times_\sk \widehat S$, where $\widehat{(\_)}$ denotes the completion of a local ring with respect to its maximal ideal. Furthermore, if $R$ and $S$ are complete, so is $P$.

\item By \cite[(1.0.3)]{AAM}, if $(A,\m_A,\sk)$ is a local ring such that $R \simeq A/I$, $S \simeq A/J$, and $\m_A = I+J$, then $R\times_{\sk} S \simeq A/I\cap J$. 

\item By \cite[Thm. 1]{DK},
\begin{equation*}\tag{\ref{FP}.1}
\frac{1}{\p^{P}(t)}=\frac{1}{\p^{R}(t)} + \frac{1}{\p^{S}(t)} -1.
\end{equation*}

\item It follows from (b) and (c) that a local ring $(P,\m_P,\sk)$ can be decomposed nontrivially as a fibre product over $\sk$ if and only if 
$\m_P =\langle y_1,\ldots,y_m, z_1, \ldots,z_n \rangle$ for $m$, $n \geq 1$ with $\langle \y \rangle \cap \langle \z \rangle = 0$. 

In this case, we see that
$P \simeq R \times_\sk S$, where $R = P/\langle \z \rangle$ and $S = P/\langle \y \rangle$.

\item Identify $\m_R$ with $\{(r,0): r \in \m_R\}$ and $\m_S$ with $\{(0,s): s \in \m_S\}$. Then $P$ is a local ring with maximal ideal $\m_P = \m_R \times \m_S$. Hence $\edim(P) = \edim(R) + \edim(S)$ and by (e), $\gr(P) \simeq \gr(R) \times_\sk \gr(S)$. Thus, for $i \geq 1$, we have $H_P(i) = H_R(i) + H_S(i)$.

\item If $R$ and $S$ are Artinian, then $\lambda(P ) = \lambda(R) + \lambda(S) - 1$. Furthermore, if $\ll(R)$, $\ll(S) \geq 1$, then $\soc(P) = \soc(R) \oplus \soc(S)$, and hence, $\type(P) = \type(R) + \type(S)$. In particular, when $R$ and $S$ are different from $\sk$, $P$ is not Gorenstein.

\end{enumerate}
}\end{remark}

\begin{lemma}[Fibre Products of Complete Regular Local Rings]\label{CP1}
Let $(\widetilde R,\m_{\widetilde R},\sk)$ and $(\widetilde S,\m_{\widetilde S},\sk)$ be complete regular local rings. Then there exists a complete regular local ring $(\widetilde P,\m_{\widetilde P},\sk)$ with $\m_{\widetilde P} = \langle Y_1, \ldots, Y_m, Z_1,\ldots, Z_n \rangle$ such that $\widetilde R \times_\sk \widetilde S \simeq {\widetilde P}/\langle \Y\cdot \Z\rangle$, $\widetilde R \simeq {\widetilde P}/\langle \Z\rangle$, and $\widetilde S \simeq {\widetilde P}/\langle \Y\rangle$. In particular, ${\widetilde P}/\langle \Y\cdot \Z\rangle$ is a Cohen presentation of $\widetilde R \times_\sk \widetilde S$.
\end{lemma}

\begin{proof}
Let $\m_{\widetilde R} = \langle y_1, \ldots, y_m \rangle$, and $\m_{\widetilde S} = \langle z_1,\ldots, z_n \rangle$. Then, with notation as in Remark \ref{FP}, the maximal ideal of $\widetilde R \times_\sk \widetilde S$ is $\m_{\widetilde R} \times \m_{\widetilde S} = \langle \y, \z \rangle$. 

    Now, $\widetilde R \times_\sk \widetilde S$ is complete by Remark \ref{FP}(b), and hence has a Cohen presentation, say ${\widetilde P}/I$. Let $\m_{\widetilde P} = \langle Y_1, \ldots, Y_m, Z_1,\ldots, Z_n \rangle$, where the $Y_i$'s and $Z_j$'s are preimages of $y_i$'s and $z_j$'s respectively. Let $J_1 = \langle \Z \rangle + I$ and $J_2 = \langle \Y \rangle + I$. Then $\widetilde R \simeq {\widetilde P}/ J_1$, $\widetilde S \simeq {\widetilde P}/ J_2$, and $\sk \simeq {\widetilde P}/(J_1 + J_2)$. Hence, $\widetilde R \times_\sk \widetilde S \simeq {\widetilde P}/(J_1 \cap J_2)$ by Remark \ref{FP}(c).

Since $\widetilde R$ is regular, the minimal generators of $\langle \Z \rangle + I$ are in $\m_{\widetilde P} \setminus \m_{\widetilde P}^2$. Hence $I \subset \m_{\widetilde P}^2$ forces $I \subset \langle \Z \rangle$. Similarly, since $\widetilde S$ is regular, $I \subset \langle \Y \rangle$. Thus $J_1 \cap J_2 = \langle \Z \rangle \cap \langle \Y \rangle = \langle \Y\cdot \Z \rangle$.
\end{proof}

The following proposition relates the defining ideals in Cohen presentations of a fibre product ring and its components.

\begin{proposition}[Cohen Presentation of a Fibre Product]\label{CP2}
Let $(R,\m_{R},\sk)$ and $(S,\m_{S},\sk)$ be complete Noetherian local rings with Cohen presentations $R \simeq \widetilde R/I_{R}$ and $S \simeq \widetilde S/I_{S}$ respectively. Let $({\widetilde P},\m_{\widetilde P},\sk)$ be as in Lemma \ref{CP1}, and $\tilde\pi_R: {\widetilde P} \rightarrow \widetilde R$ and $\tilde\pi_S: {\widetilde P} \rightarrow \widetilde S$ be the natural projections.

If $P \simeq R \times_\sk S$, $J_R = \tilde\pi_R^{-1}(I_{R})$, and $J_S = \tilde\pi_S^{-1}(I_{S})$,
then $P = {\widetilde P}/I_P$ where
\begin{enumerate}[{\rm a)}]
\item $\langle \Y \cdot \Z \rangle \subset I_P \subset \langle \Y, \Z \rangle^2$. In particular, ${\widetilde P}/I_P$ is a Cohen presentation for $P$.

\item $I_P = (J_R \cap \langle \Y \rangle)  + (J_S \cap \langle \Z \rangle) + \langle \Y\cdot\Z \rangle$.

\item $J_{R}= (I_P \cap \langle \Y \rangle) + \langle \Z \rangle$, and $J_{S} = (I_P \cap \langle \Z \rangle) + \langle \Y \rangle$.

\item $\mu(I_P) = \mu(I_{R}) + \mu(I_{S}) + mn$.
\end{enumerate}
\end{proposition}

\begin{proof}
Let the notation be as in the previous lemma. We have the following commutative diagram:
 \begin{equation}
\begin{gathered}
\xymatrixrowsep{1pc}
\xymatrixcolsep{0.9pc}
\xymatrix{
&&& \widetilde R
\ar@{->}[rr]
&& R
\ar@{->}[drr]
 \\
{\widetilde P} 
\ar@/^1pc/[urrr]^-{\tilde\pi_R}
\ar@{->}[rr]
\ar@/_1pc/[drrr]_-{\tilde\pi_S}
&& \widetilde R\times_k\widetilde S
\ar@{->}[ur]
\ar@{->}[rr]
\ar@{->}[dr]
&& P = R\times_kS
\ar@{->}[ur]
\ar@{->}[rrr]
\ar@{->}[dr]
&& &k
 \\
&&& \widetilde S
\ar@{->}[rr]
&& S
\ar@{->}[urr]
}
  \end{gathered}
\end{equation}
where all the arrows denote canonical
surjections.

Then $R \simeq {\widetilde P}/J_R$, $S \simeq {\widetilde P}/J_S$, and since $\langle \Z \rangle \subset J_R$ and $\langle \Y \rangle \subset J_S$, we have $\sk \simeq {\widetilde P}/(J_R + J_S)$. Hence, by Remark \ref{FP}(c), we can write $P \simeq {\widetilde P}/I_P$, where $I_P = J_R \cap J_S$. 

(a) By Lemma \ref{CP1} and Remark \ref{FP}(f), $\edim({\widetilde P}) = \edim({\widetilde R}\times_\sk {\widetilde S}) = \edim(P)$, which forces $I_P \subset \m_{\widetilde P}^2$. Furthermore, since the map from ${\widetilde P}$ to $P$ factors through $\widetilde R\times_k\widetilde S$, we have $\langle \Y\cdot \Z\rangle \subset I_P$. 

(b) Since $\langle \Y \rangle \subset J_S$ and $\langle \Z \rangle \subset J_R$, we see that $J_R \cap \langle \Y \rangle \subset J_R \cap J_S = I_P$. Similarly, $J_S \cap \langle \Z \rangle \subset I_P$. This proves one inclusion in (b).

Now, let $F \in I_P \subset \m_{\widetilde P}^2$, $F_Y \in \langle \Y \rangle^2$ and $F_Z \in \langle \Z \rangle^2$ be such that $F - F_Y - F_Z \in \langle \Y\cdot \Z\rangle \subset I_P$.
Since $F_Z \in \langle \Z \rangle \subset J_R$, $F - F_Y - F_Z$ and $F \in I_P \subset J_R$, we have $F_Y \in J_R$. Thus $F_Y \in J_R \cap \langle \Y \rangle$. Similarly, $F_Z \in J_S \cap \langle \Z \rangle$, proving (b).

(c) Note that $J_R = I_P + \langle \Z \rangle$. Hence, we have $(I_P \cap \langle \Y \rangle) + \langle \Z \rangle \subset J_{R}$. Now, let $F \in J_{R}$, and $F_Z \in \langle \Z \rangle$ be such that $F - F_Z = G \in I_P$. Write $G = G_Y + G_Z$, where $G_Y \in \langle \Y \rangle$ and $G_Z \in \langle \Z \rangle$. Now, $F$ and $F_Z + G_Z \in J_{R}$ force $G_Y \in J_{R}$. 
Since $G_Y \in \langle \Y \rangle \subset I_P$, we get $J_{R} \subset (I_P \cap \langle \Y \rangle) + \langle \Z \rangle$. 
Thus $J_{R}= (I_P \cap \langle \Y \rangle) + \langle \Z \rangle$, and by symmetry, $J_{S}= (I_P \cap \langle \Z \rangle) + \langle \Y \rangle$.

(d) Comparing the coefficients of $t$ and $t^2$ in Equation (\ref{FP}.1), we see that $\beta_1^P = \beta_1^{R} + \beta_1^{S}$
and $\beta_2^P = \beta_2^{R} + \beta_2^{S} + 2 \beta_1^{R}\beta_1^{S}$. Hence, by Remark \ref{PS2}, we see that $\mu(I_P) = \mu(I_{R}) + \mu(I_{S}) + mn$.
\end{proof}

\begin{remark}\label{CP3}{\rm
From their definitions, $J_R = (J_R \cap \langle \Y \rangle) + \langle \Z \rangle$, and $J_S = (J_S \cap \langle \Z \rangle) + \langle \Y \rangle$. This observation plays an important role in the proof of Proposition \ref{Prop2}.
}\end{remark}

\subsection{Connected Sums}

As we see in Remark \ref{FP}(g), if $(R,\m_R,\sk)$ and $(S,\m_S,\sk)$ are Artinian local rings, neither of which is a field, then $P = R \times_\sk S$ is not Gorenstein. We define an appropriate quotient called a connected sum which is Gorenstein. More details can be found in \cite[Section 2]{AAM} and \cite[Chapter 4]{HAthesis}.

\begin{definition}\label{CSdef}{\rm 
Let $(R,\m_R,\sk)$ and $(S,\m_S,\sk)$ be Gorenstein Artin local rings different from $\sk$.
Let $\soc(R) = \langle\delta_R\rangle$, $\soc(S) = \langle\delta_S\rangle$. Identifying $\delta_R$ with
$(\delta_R, 0)$ and $\delta_S$ with $(0, \delta_S)$, a {\it connected sum} of $R$ and $S$ over $\sk$,
denoted $R\#_\sk S$, is the ring $R\#_\sk S = (R\times_\sk S)/\langle\delta_R - \delta_S\rangle$.
}\end{definition}

Connected sums of $R$ and $S$ over $\sk$ depend on
the generators of the socle $\delta_R$ and $\delta_S$ chosen. For example, the connected sums 
$Q_1 = (R \times_\sk S) / \langle y^2 - z^2\rangle$ and $Q_2 = (R \times_\sk S) /\langle y^2 - 5 z^2\rangle$ of 
$R = \QQ[Y]/\langle Y^3\rangle$ and $S = \QQ[Z]/\langle Z^3\rangle$ are not isomorphic as rings, as shown in \cite[Ex. 3.1]{AAM}.

\begin{remark}\label{CS} {\rm With notation as in Definition \ref{CSdef}, set $P = R \times_{\sk} S$ and let
$Q = R\#_\sk S $.

\begin{enumerate}[a)]
\item\label{trivial} (Trivial connected sum). Every Gorenstein Artin local ring, which is not a field, is trivially a connected sum over its residue field. 
In order to see this, consider a Gorenstein Artin local ring $(R,\m_R,\sk)$ with $\ll(R) \geq 1$ and let
$S$ be a $\sk$-algebra of length two. Note that this forces $S$ to be Gorenstein.
One can check that $R \#_\sk S \simeq R$.

\item Note that $\lambda(Q) = \lambda(P) - 1 = \lambda(R) + \lambda(S) - 2$ since $0 \neq \delta_R - \delta_S \in \soc(P)$.

\item If $\ll(R)$, $\ll(S) \geq 2$, then
$\edim(Q) = \edim(R) + \edim(S)$.

\item If $\ll(R)$, $\ll(S) \geq 1$, then $Q$ is a Gorenstein Artin local ring with $\ll(Q) = \max\{\ll(R),\ll(S)\}$. 
This is proved in \cite[Prop. 4.4]{Le}; see \cite[Thm. 2.8]{AAM} for a more general result.

\item By the definition of $Q$, it is clear that $\overline{Q} \simeq \overline P \simeq \overline R \times_\sk \overline S$, where $\ \bar{ }\ $ 
denotes going modulo the respective socles. We prove a partial converse of this observation in Proposition \ref{Socle}.

\item If $\edim(R) = m$, $\edim(S) = n$, $\ll(R)$, $\ll(S) \geq 2$, then by (d), (e) and Remark \ref{PS}, we get
\begin{equation*}\tag{\ref{CS}.1}
\frac{1}{\p^{Q}(t)}=\frac{1}{\p^{R}(t)} + \frac{1}{\p^{S}(t)} - 1 + \phi(m,n)t^2 = \frac{1}{\p^P(t)} + \phi(m,n)t^2,
\end{equation*}
where $\phi(m,n) = - 1$ when $m$, $n \geq 2$, $\phi(1,1) = 1$,
and $\phi(m,n) = 0$ otherwise.

\item\label{CS1}(Cohen presentation).
Let $R \simeq \widetilde P/J_R$, $S \simeq \widetilde P/J_S$ and $P \simeq \widetilde P/I_P$ be as in Proposition \ref{CP2}.

Suppose $R$ and $S$ are Gorenstein and $Q = R \#_\sk S$ is their connected sum over $\sk$. Then, by
Definition \ref{CSdef}, since $\langle \Z \rangle \subset J_R$ and $\langle \Y \rangle \subset J_S$, there exists $\Delta_R \in \langle\Y\rangle$ and 
$\Delta_S \in \langle\Z\rangle$ such that their respective
images $\delta_R \in R$ and $\delta_S \in S$ generate the respective socles and
$Q \simeq  P/\langle\delta_R - \delta_S\rangle$.
Thus $Q \simeq \widetilde P/I_Q$, where

$I_Q= I_P + \langle \Delta_R-  \Delta_S\rangle =
(J_R \cap \langle \Y \rangle)  + (J_S \cap \langle \Z \rangle) + \langle \Y\cdot\Z \rangle+ \langle \Delta_R-  \Delta_S\rangle$.

In particular, if $R$ and $S$ are standard graded $\sk$-algebras, then $Q$ is standard graded if and only if $R$ and $S$ have the same Loewy length. In this case, $\ll(Q) = \ll(R) = \ll(S)$.

\item By (g), if $Q$ is a connected sum, then we can write $\m_Q = \langle \y, \z \rangle$, with $\y \cdot \z = 0$. Furthermore, since $Q$ is Gorenstein, and hence not decomposable as a fibre product, $\langle \y \rangle \cap \langle \z \rangle = \soc(Q)$.
\end{enumerate}
}\end{remark}

\section{Connected Sums and Indecomposibility}\label{Indec}

The main question we would like to address is:

\begin{question}[Main Question]\label{MQ}{\rm
When is a Gorenstein Artin local ring decomposable as a connected sum over its residue field?
}\end{question}

In light of Remark \ref{CS}(\ref{trivial}), we make the following key definition about connected sums. A similar terminology is also used
for fibre products in this article.

\begin{definition}{\rm
Let $(Q,\m,\sk)$ be a Gorenstein Artin local ring. We say that $Q$ {\it decomposes as a connected sum}
over $\sk$ if there exist Gorenstein Artin local rings $R$ and $S$ such that
$Q \simeq R \#_\sk S$ and $R \not\simeq Q \not\simeq S$.
In this case, we call $R$ and $S$ the {\it components} in a connected sum decomposition of $Q$, and say that
$Q \simeq R \#_\sk S$ is a non-trivial decomposition.

If $Q$ cannot be decomposed as a connected sum over $\sk$,
we say that $Q$ is {\it indecomposable} as a connected sum over $\sk$.
}\end{definition}

\begin{remark}{\rm
For Gorenstein Artin rings $(R,\m_R,\sk)$ and $(S,\m_S,\sk)$ with $\ll(R)$, $\ll(S) \geq 2$, we see that $Q = R \#_\sk S$ is
a non-trivial decomposition of $Q$ as a connected sum over $\sk$.

Indeed, $\lm(Q) = \lm(R) + \lm(S) -2$ and $\ll(R)$, $\ll(S) \geq 2$
force $\lm(Q) > \max\{\lm(R), \lm(S)\}$, hence $R \not\simeq Q \not\simeq S$.
}\end{remark}

Let $(Q,\m,\sk)$ be a Gorenstein local ring. It is known (see \cite[8.3]{AAM}) that if $Q$ is a {\it complete intersection} with $\edim(Q) \geq 3$, then $Q$ is
indecomposable as a connected sum over $\sk$. We give an easier proof in the Artinian case, which follows from Proposition \ref{mingensProp}, as is noted in Theorem \ref{mingensTheorem}.

In this section, we see a condition on  the Hilbert function of $Q$ (Theorem \ref{HilbertCoeff2}), which forces indecomposibility.
When the embedding dimension of $Q$ is three or four, one can give
conditions in terms of the minimal number of generators of the defining ideal of $Q$ which force indecomposibility
(Theorem \ref{mingensTheorem}).

The following remark is a quick detour relating the notion of connected sums and indecomposibility of $\sk$-algebras with {\it Macaulay's inverse systems}.

\begin{remark}[Connected Sums and Inverse Systems]\label{apolar}{\rm
 A different point of view to study Gorenstein Artin $\sk$-algebras is via inverse systems,
in which such rings correspond to polynomials. For more details, see \cite[Section 1.4]{HAthesis} or \cite[Section 2]{ER}.

If $R$ and $S$ are Gorenstein Artin $\sk$-algebras corresponding to polynomials $F(\Y)$ and
$G(\Z)$ respectively, then the Gorenstein Artin $\sk$-algebra corresponding to $F + G$ is a connected sum of $R$ and
$S$ over $\sk$. (For example, see \cite[4.24]{HAthesis}).

In terms of inverse systems, answering the Question \ref{MQ} amounts to the following:
Given a polynomial $F$ corresponding to $Q$, write $F = F_1 + F_2$, where
$F_1$ and $F_2$ are polynomials in disjoint sets of variables.

When $F$ is homogeneous, the above property has been studied in \cite{BBKT}. The authors define such a polynomial
to be a {\it direct sum} and the corresponding Gorenstein algebra to be {\it apolar}.
}\end{remark}

\subsection{Minimal Number of Generators of the Defining Ideal}\hfill{}

We begin our study of indecomposibility with a result on the minimal number of generators of the defining ideal, $I_Q$, of $Q$.

\begin{proposition}\label{mingensProp} Let $Q$, $R$ and $S$ be Gorenstein Artin local rings with the same residue field $\sk$, and respective Cohen presentations $\widetilde Q/I_Q$, 
$\widetilde R/I_R$, and $\widetilde S/I_S$.
If $Q \simeq R \#_\sk S$ is a non-trivial decomposition, then we have
$$\mu(I_Q) = \mu(I_R) + \mu(I_S) + mn + \phi(m,n),$$

where $\phi(m,n) = 1$ when $m$, $n \geq 2$, $\phi(1,1) = -1$,
and $\phi(m,n) = 0$ otherwise.
\end{proposition}

\begin{proof}
First of all, note that if $m = 1 = n$, then $\mu(I_R)
= 1 = \mu(I_S)$. Since $Q$ is Gorenstein Artin with $\edim(Q) = m +
n = 2$, a well known result of Serre (see \cite{Se}) shows that $\mu(I_Q) = 2$. Hence, without loss of
generality, we may assume that $m \geq 2$.

It is clear from Equation \ref{CS}.1 that $\beta_2^Q = \beta_2^P$ when $n = 1$ and $\beta_2^Q = \beta_2^P + 1$
when $n \geq 2$. Since $\beta_1^Q = \beta_1^P$ in either case, the proof is complete by Remark \ref{PS2}.
\end{proof}

As a consequence, we can prove:

\begin{theorem}\label{mingensTheorem} A Gorenstein Artin local ring $(Q,\m_Q,\sk)$, with Cohen presentation $\widetilde Q/I_Q$, is indecomposable as a
connected sum over $\sk$ when one of the following holds:

\begin{enumerate}[{\rm a)}]
\item $\edim(Q) = 3$ and $\mu(I_Q) \neq 5$.

\item $\edim(Q) = 4$ and $\mu(I_Q)$ is an even number.

\item $\edim(Q) \geq 3$ and $Q$ is a complete intersection ring.
\end{enumerate}
\end{theorem}

\begin{proof}
Suppose $Q \simeq R \#_\sk S$ is a non-trivial decomposition of $Q$ as a connected sum over $\sk$, with $\edim(R) = m$ and $\edim(S) = n$. Then $d = \edim(Q) = m+n$.

(i) Suppose $\edim(Q) = 3$. Then without loss of generality, $m = 2$ and $n = 1$. Hence
$\mu(I_S) = 1$ and by the above-mentioned result of Serre, $\mu(I_R) = 2$.
Thus $\mu(I_Q) = 5$.

(ii) If $\edim(Q) = 4$, then either $m = n = 2$, in which case $\mu(I_R) = 2 = \mu(I_S)$ forcing $\mu(I_Q) = 9$
or, without loss of generality, $m = 3$ and $n=1$. In this case, $\mu(I_S) = 1$, and by \cite{Wa}, $\mu(I_R)$ is an odd number.
Hence, in either case, by Proposition \ref{mingensProp}, $\mu(I_Q)$ is odd.

(iii) Since $\mu(I_R) \geq m$, $\mu(I_S) \geq n$, and $mn \geq 1$, we observe that $\mu(I_Q) = m + n$ if and only if $\mu(I_R) = m$, $\mu(I_S) = n$ and $mn = 1$, i.e., when $m = 1$, $n = 1$ and $\edim(Q) = 2$.
\end{proof}

\begin{example}{\rm
Let $R \simeq \QQ[Y]/\langle Y^3 \rangle$ and $S \simeq \QQ[Z]/\langle Z^3 \rangle$. Then $Q = R \#_\QQ S$ is decomposable as a connected sum, but is a complete intersection ring, since $\edim(Q) = 2$. This shows that the condition $\edim(Q) \geq 3$ is necessary in Theorem \ref{mingensTheorem}(iii).
}\end{example}

The above theorem does not give necessary conditions, see Example \ref{mingensExample}.

\subsection{Hilbert Functions}
We first obtain a numerical criterion satisfied by connected sums.

\begin{proposition}\label{HilbertCoeff1}
Let  $Q = R\#_\sk S$, where $(R,\m_R,\sk)$ and $(S,\m_S,\sk)$ are Gorenstein Artin local rings, with $\ll(R)$, $\ll(S) \geq 2$, $\edim(R) = m$ and $\edim(S) = n$.
Then $H_Q(2) \leq \binom{m+n+1}{2} - mn$.
\end{proposition}

\begin{proof}
With notation as in Remark \ref{CS}(\ref{CS1}), since $\Y \cdot \Z \subset I_Q$, we have
$\Y^* \cdot  \Z^* \subset I_Q^*$. Hence, $H_{\widetilde Q}(2) = \binom{m+n+1}{2}$ implies that
\vskip 2pt
$H_Q(2) =  H_{\gr(Q)}(2) = \binom{m+n+1}{2} -
\dim_{\sk}((I_Q^* + \langle\Y,\Z\rangle^3)/ \langle\Y, \Z\rangle^3)
\leq \binom{m+n+1}{2} - mn$.
\end{proof}

 This gives us a sufficient condition for indecomposibility in terms of Hilbert functions:

\begin{theorem}\label{HilbertCoeff2}
Let $(Q,\m_Q,\sk)$ be a Gorenstein Artin local ring with $\edim(Q) = d$. If for some $i \geq 2$, $H_Q(i) \geq \binom{d - 2 + i}{i} + 2$, then $Q$ is
indecomposable as a connected sum over $\sk$.
\end{theorem}

\begin{proof}
Note that for positive integers $m$ and $n$, if $d = m+n$ is fixed, the minimum value of $mn = m(d-m)$ is obtained when
$m = 1$ or $m = d-1$. Hence,  by Proposition \ref{HilbertCoeff1}, to prove the indecomposibility of $Q$, it is enough to show that
$H_Q(2)  > \binom{d+1}{2} - (d -1) = \binom{d}{2} + 1$.

If $H_Q(2) \leq \binom{d}{2} + 1$, then by Macaulay's theorem on Hilbert functions (e.g., see \cite[6.3.8]{HH}), and induction on $i$, it can be seen that
$H_Q(i) \leq \binom{d - 2 + i}{i} + 1$ for each $i \geq 3$. Thus, if $H_Q(i) \geq \binom{d - 2 + i}{i} + 2$ for some $i \geq 2$, then $H_Q(2) \geq \binom{d}{2} + 2$, and hence
$Q$ is indecomposable as a connected sum over $\sk$.
\end{proof}

Theorem \ref{mingensTheorem}, and the condition on $H_Q(i)$ in Theorem \ref{HilbertCoeff2}, do not give necessary conditions for indecomposibility, as can be seen from the  following examples.

\begin{example}\label{mingensExample}\hfill{}{\rm

1) Let $Q \simeq \QQ[Y_1,Y_2,Y_3]/I_Q$, where $I_Q = \langle Y_1^4, Y_2^4, Y_3^4,Y_1^2Y_2^2 - Y_1^2Y_3^2, Y_1^2 Y_2^2 - Y_2^2Y_3^2 \rangle$. Then $Q$ is Gorenstein with $d = \edim(Q) = 3$. Furthermore, $H_Q(2) = 6 \geq \binom{d}{2} + 2$, hence by Theorem \ref{HilbertCoeff2}, $Q$ is indecomposable as connected sum over $\QQ$.
However, $\mu(I_Q) = 5$. 

2) Let $Q = \QQ[Y_1,Y_2,Y_3,Y_4]/I_Q$, where $I_Q = \langle Y_1^4,Y_2^4,Y_3^4,Y_4^4, Y_1^2Y_2^2 - Y_i^2Y_j^2: 1 \leq i < j \leq 4 \rangle$. Then $Q$ is Gorenstein,
$d = \edim(Q) = 4$, $\mu(I_Q) = 9$ is odd, and $Q$ is indecomposable as a connected sum since $H_Q(2) = 10 \geq \binom{d}{2} + 2$.

3) Let $Q = \QQ[X_1,X_2]/\langle X_1^2X_2, X_1^3 - X_2^2\rangle$. Then $Q$ is Gorenstein, $H_Q(2) = 2 < \binom{2}{2} +2$, and $\edim(Q) = 2$. However, $Q$ is indecomposable as a connected sum over $\QQ$. 

Indeed, if $Q \simeq R \#_\QQ S$ is a non-trivial decomposition, for some Gorenstein Artin $\QQ$-algebras $R$ and $S$, then, by Remark \ref{CS}(\ref{CS1}), one can find elements $Y$ and $Z$ in
$\langle X_1, X_2 \rangle \setminus \langle X_1, X_2 \rangle^2$ such that $Y \cdot Z \in I_Q$ and $\langle Y, Z \rangle = \langle X_1, X_2 \rangle$.

Now, $Y \cdot Z \in I_Q=\langle X_1^2X_2, X_1^3 - X_2^2\rangle$ implies that $Y Z + c X_2^2 \in \langle X_1, X_2 \rangle^3$ for some $c \in \QQ$. We write $Y = a_1X_1 + a_2 X_2 + F$ and $Z = b_1X_1 + b_2 X_2 + G$, where $a_1$, $a_2$, $b_1$, $b_2 \in \QQ$, and $F$, $G \in \langle X_1, X_2 \rangle^2$.
Then $Y Z + c X_2^2 \in  \langle X_1, X_2 \rangle^3$ forces $(a_1,  b_1)=(0,0)$, $(a_1, a_2)=(0,0)$, or $(b_1, b_2)=(0,0)$  contradicting $\langle Y, Z \rangle = \langle X_1, X_2 \rangle$.
}\end{example}

A Gorenstein Artin local ring $Q$ is said to be {\it compressed} if it has a maximum possible Hilbert function given the
embedding dimension $d$ and Loewy length $s$, i.e., if the Hilbert function of $Q$ is
$H_Q(i)=\text{min}\{\binom{d+i-1}{i}, \binom{d+s-i-1}{s-i}\}$.

\begin{corollary}\label{compressed}
If $(Q,\m_Q,\sk)$ is a compressed Gorenstein Artin local ring with $\ll(Q) \geq 4$, then $Q$ is
indecomposable as a connected sum over $\sk$.
\end{corollary}

\begin{remark}\label{generic}{\rm
Since {\it generic Gorenstein $\sk$-algebras} are compressed (see \cite[Thm. 1]{Ia84}), they are
indecomposable as connected sums over $\sk$.
This was proved in \cite[4.4]{SS} by using different techniques.
}\end{remark}

The following example shows that there are rings which are neither complete intersections, nor compressed algebras, which are indecomposable as a connected sum.
\begin{example}
Let $\sk = \ZZ/2 \ZZ$, and $Q = \sk[X_1,X_2,X_3]/\langle X_1^3,X_2^3,X_3^3,X_1X_2X_3,X_1^2+X_2^2+X_3^2 \rangle$. Then $Q$ is Gorenstein, which is clearly not a complete intersection. Furthermore, $\ll(Q) = 4$, and $H_Q(2) = 5$. Thus, $H_Q(2) < 6$ implies $Q$ is not compressed, and $H_Q(2) > 4$ forces $Q$ to be indecomposable as a connected sum over $\sk$.
\end{example}

\begin{remark}[Almost Stretched Rings]{\rm
It follows from Sally's work (\cite{Sa1}) that a {stretched} Gorenstein Artin ring, of embedding dimension at least 2, is decomposable as a connected sum. A Gorenstein Artin local ring $(Q,\m_Q,\sk)$ is {\it almost stretched} if $\mu(\m_Q^2) \leq 2$. 

In a private conversation, Paolo Mantero asked if the same is true in the {almost stretched} case. In Example \ref{mingensExample}(3), $Q$ is almost stretched, but is indecomposable, answering the question in the negative.
}\end{remark}

\section{Criteria for Decomposability} \label{Dec}
\subsection{A criterion in terms of the defining ideals}

The following gives a condition on the defining ideal of a Gorenstein Artin local ring which forces it to be a connected sum.

\begin{proposition}\label{MainProp}
Let $(Q,\m_Q,\sk)$ be a Gorenstein Artin local ring and $Q = \widetilde Q/I_Q$ be its Cohen presentation, with $\m_{\widetilde Q} = \langle Y_1, \ldots, Y_m, Z_1, \ldots, Z_n\rangle$, $m$, $n \geq 1$. Let $R = \widetilde Q/J_R$ and
$S = \widetilde Q/J_S$, where $J_R = (I_Q \cap \langle \Y \rangle) + \langle \Z \rangle$ and $J_S = (I_Q \cap \langle\Z\rangle) + \langle \Y \rangle$. Suppose
$\Y \cdot \Z \subset I_Q$. Then
\begin{enumerate}[{\rm a)}]
\item $I_Q \cap \langle \Y \rangle = J_R \cap \langle \Y \rangle$ and $I_Q \cap \langle \Z \rangle = J_S \cap \langle \Z \rangle$.
\item $R$ and $S$ are Gorenstein Artin and
\item $Q \simeq R \#_\sk S$.
\end{enumerate}
\end{proposition}

\begin{proof}
(a) Clearly $I_Q \cap \langle \Y \rangle \subset J_R \cap \langle \Y \rangle$. Now, let $\Lambda = \Lambda_Y + \Lambda_Z \in J_R \cap \langle \Y \rangle$, where $\Lambda_Y \in I_Q \cap \langle \Y \rangle$ and $\Lambda_Z \in \langle \Z \rangle$. Since $\Lambda_Z = \Lambda - \Lambda_Y \in \langle \Y \rangle \cap \langle \Z \rangle = \langle \Y \cdot \Z \rangle \subset I_Q \cap \langle \Y \rangle$, we see that $\Lambda \in I_Q \cap \langle \Y \rangle$. Thus $I_Q \cap \langle \Y \rangle = J_R \cap \langle \Y \rangle$. Similarly, $I_Q \cap \langle \Z \rangle = J_S \cap \langle \Z \rangle$.

(b) By symmetry, it suffices to prove (b) for $R$. Let $0 \neq \delta_R \in \soc(R)$, and $\Delta_R \in \langle \Y \rangle \setminus J_R$ be a preimage of $\delta_R$ in $\widetilde Q$. We want to prove that $\langle \Delta_R \rangle R = \soc(R)$. Since $\Delta_R \not\in J_R$, and $\Delta_R \in \langle \Y \rangle$, we see, by (a), that $\Delta_R \not\in I_Q$.

Now  
$\Y \cdot \Delta_R \subset J_R \cap \langle \Y \rangle \subset I_Q$. Moreover, $\Y \cdot \Z \subset I_Q$ implies 
$\Z \cdot \Delta_R \subset I_Q$. Hence $\langle \Delta_R\rangle Q \subset \soc(Q)$. Since $\Delta_R \not\in I_Q$, and $\dim_{\sk}(\soc(Q)) = 1$, we see that $\langle \Delta_R \rangle Q = \soc(Q)$. 

Let $\delta \in \soc(R)$ and $\Delta \in \langle \Y \rangle$ be a preimage in $\widetilde Q$. As seen above, $\langle \Delta \rangle Q \subset \soc(Q) = \langle \Delta_R \rangle Q$, i.e., $\Delta \in \langle \Delta_R \rangle + I_Q$. Since $\Delta$, $\Delta_R \in \langle \Y \rangle$, we have $\Delta \in \langle \Delta_R \rangle + (\langle \Y \rangle \cap I_Q)  \subset  \langle \Delta_R \rangle + J_R$ by (a). Thus $\delta \in \langle \delta_R \rangle$, which implies that $\lambda(\soc(R)) = 1$, proving that $R$ is Gorenstein.

(c) Let $P = R \times_\sk S$. Then $P \simeq \widetilde Q/I_P$, where
$I_P = (J_R \cap \langle \Y \rangle) + (J_S \cap \langle \Z \rangle) + \langle \Y \cdot \Z \rangle$. By the hypothesis and (a), $I_P \subset I_Q$ and hence
there is a natural surjective map $\pi: P \longrightarrow Q$.

Let $\Delta_R \in \langle \Y \rangle$ and $\Delta_S \in \langle \Z \rangle$ be such that $\soc(R) = \langle \Delta_R \rangle R$ and $\soc(S) = \langle \Delta_S \rangle S$. Let their corresponding images in $P$ be $\rho_R$ and $\rho_S$.
The same approach as in (b) shows that $\langle \Delta_R, \Delta_S \rangle P = \langle \rho_R,\rho_S \rangle \subset \soc(P)$.
We claim that:
\smallskip

$Claim. \quad \quad \quad 0 \subsetneq \langle \rho_R - u \rho_S\rangle \subset \ker(\pi) \subsetneq \langle \rho_R, \rho_S\rangle = \soc(P).$

\smallskip
We first show that $\soc(P) = \langle \rho_R, \rho_S \rangle$. To see this, let $\Delta \in \widetilde Q$ be such that $\langle \Y, \Z \rangle \Delta \subset I_P$. Write $\Delta = \Delta_Y + \Delta_Z$, where $\Delta_Y \in \langle \Y \rangle$ and $\Delta_Z \in \langle \Z \rangle$, with images $\rho$, $\rho_Y$ and $\rho_Z$ in $P$ respectively.

Observe that $\Y \cdot \Z \subset I_P$ implies that $\Y \Delta_Y \in I_P \cap \langle \Y \rangle \subset J_R$. Since $\langle \Z \rangle \subset J_R$, we have $\langle \Delta_Y \rangle R \subset \soc(R) = \langle \Delta_R \rangle R$. In particular, $\Delta_Y \in \langle \Delta_R \rangle + (J_R \cap \langle \Y \rangle) \subset \langle \Delta_R \rangle + I_P$, and hence $\rho_Y \in \langle \rho_R \rangle$ in $P$. A similar proof shows that $\rho_Z \in \langle \rho_S \rangle$ in $P$, proving the last equality in the claim. 

In order to prove that $\ker(\pi) \subset \soc(P)$, let $\gamma \in \ker(\pi)$ and $\Gamma \in I_Q$ be a preimage in 
$\widetilde Q$. We need to show that $Y_i \Gamma$, $Z_j \Gamma \in I_P$ for all $i$ and $j$. To do this, it is enough to prove that $Y_1 \Gamma \in I_P$. Write 
$\Gamma = \Gamma_Y + \Gamma_Z$, where $\Gamma_Y \in \langle \Y \rangle$ and $\Gamma_Z \in \langle \Z \rangle$. Then $Y_1 \Gamma_Z \in I_P \subset I_Q$.

Thus $Y_1 \Gamma_Y = Y_1 \Gamma - Y_1\Gamma_Z \in I_Q \cap \langle \Y\rangle = J_R \cap \langle \Y \rangle \subset I_P$. Therefore,
$Y_1 \Gamma \in I_P$ proving that $\ker(\pi) \subset \soc(P)$. Observe that $\Delta_R \not\in I_Q$ implies that $\rho_R \not\in \ker(\pi)$, and similarly, $\rho_S \not\in \ker(\pi)$. In particular, $\ker(\pi) \subsetneq \soc(P)$.

Now, as observed in the proof of (b), $\soc(Q) = \langle \Delta_R \rangle Q$, and similarly, $\soc(Q) = \langle \Delta_S \rangle Q$, i.e., $\Delta_R - U \Delta_S \in I_Q$ for some $U \in \widetilde Q$. Thus, if $u$ is the image of $U$ in $P$, then $\rho_R - u \rho_S \in \ker(\pi)$. Furthermore, $\rho_R \not\in \ker(\pi)$ forces $u$ to be a unit in $P$.

Finally, by Remark \ref{FP}(g), $\soc(P) = \langle \rho_R,\rho_S \rangle$ is not cyclic, and hence $\rho_R - u \rho_S \neq 0$ in $P$. This completes the proof of the claim.

Since $\lambda(\soc(P)) = 2$, the claim forces $\ker(\pi) = \langle \rho_R - u \rho_S\rangle$. Therefore, 
$Q \simeq P/\langle \rho_R - u \rho_S\rangle$ for some unit $u \in P$, where $\langle \rho_R \rangle R = \soc(R)$, and $\langle u\rho_S \rangle S = \soc(S)$. Hence, $Q$ is a connected sum of 
$R$ and $S$ over $\sk$, by Definition \ref{CSdef}.
\end{proof}

If $(Q,\m_Q,\sk)$ is a Gorenstein Artin local ring, and $R$ and $S$ are defined as in the above proposition,  they may not be Gorenstein in general. The requirement that $\Y \cdot \Z \subset I_Q$ is necessary for condition (b) in the above proposition to be true, as can be seen from the following example.

\begin{example}{\rm
Let $Q=\QQ[|Y,Z_1,Z_2|]/I_Q$, where $I_Q=\langle Y^2-Z_1Z_2, Z_1^3, Z_2^2\rangle$. Then $Q$ is Gorenstein.
Note that $S = \QQ[|Y,Z_1,Z_2|]/((I_Q \cap \langle Z_1,Z_2\rangle) + \langle Y \rangle) \simeq \QQ[|Z_1,Z_2|]/\langle Z_1^3,Z_2^3,Z_1^2Z_2^2 \rangle$.
Thus $S$ is not Gorenstein because $\soc(S) = \langle Z_1^2Z_2,Z_1Z_2^2\rangle$ is a two-dimensional vector space over $\sk$.
}\end{example}

The following theorem gives an analogue of Remark \ref{FP}(e) for decomposition as a connected sum over $\sk$. The proof follows from Remark \ref{CS}(g), and Proposition \ref{MainProp}.

\begin{theorem}\label{MainThm}
Let $(Q,\m_Q,\sk)$ be a Gorenstein Artin local ring. Then $Q$ can be decomposed nontrivially as a connected sum over $\sk$ if and only if $\m_Q = \langle y_1, \ldots, y_m, z_1,\ldots, z_n \rangle$, $m, n \geq 1$, such that $\y\cdot \z  = 0$.
\end{theorem}

By definition of the fibre product, if $P = R\times_\sk S$, then $R$ and $S$ can be identified with quotients of $P$. On the other hand,
if $Q = R \#_\sk S$, it is not clear how one can recover the components $R$ and $S$ from
$Q$. In the Artinian case, Proposition \ref{Prop2} allows one to find defining ideals for $R$ and $S$ in terms of a Cohen presentation of $Q$. 

\begin{proposition} \label{Prop2}
Let $R$ and $S$ be Gorenstein Artin local rings with $Q = R \#_\sk S$. Then there is a regular local ring $(\widetilde Q, \m_{\widetilde Q},\sk)$ with $\m_{\widetilde Q} = \langle \Y, \Z \rangle$, such that $Q \simeq \widetilde Q/I_Q$, $R \simeq \widetilde Q /J_R$, and $S \simeq \widetilde Q/J_S$, where $\langle \Y \cdot \Z\rangle \subset I_Q \subset \langle \Y , \Z\rangle^2$, $J_R = (I_Q\cap \langle \Y \rangle)+ \langle \Z\rangle$, and $J_S = (I_Q\cap \langle \Z \rangle)+ \langle \Y\rangle$.
\end{proposition}

\begin{proof} Firstly note that by Remark \ref{CS}(\ref{CS1}), we have a Cohen presentation $\widetilde Q/I_Q$ of $Q$ such that $\m_{\widetilde Q} = \langle \Y, \Z \rangle$, 
$R \simeq \widetilde Q /J_R$, and $S \simeq \widetilde Q/J_S$, with $(J_R \cap \langle \Y \rangle) + (J_S \cap \langle \Z \rangle) + \langle \Y \cdot \Z\rangle \subset I_Q$.

Now, let $J_{R'} = (I_Q \cap \langle \Y\rangle)+ \langle \Z \rangle$ and $J_{S'} = (I_Q \cap \langle \Z\rangle) + \langle \Y \rangle$. The proof is complete if we prove $J_R = J_{R'}$, and $J_S = J_{S'}$.

By Remark \ref{CP3}, we have $J_R \subset J_{R'}$, and $J_S \subset J_{S'}$, which induce natural surjective maps
$\pi_1: R \longrightarrow R'$ and $\pi_2: S \longrightarrow S'$, where $R' = \widetilde Q/J_{R'}$, and $S' = \widetilde Q/I_{S'}$. 
In particular, $\lambda(R') \leq \lambda(R)$ and $\lambda(S') \leq \lambda(S)$. In order to prove that
$J_R = J_{R'}$ and $J_S = J_{S'}$, it is enough to prove that $\pi_1$ and $\pi_2$ are isomorphisms, in particular, it is enough to
show that $\lambda(R') = \lambda(R)$ and $\lambda(S') = \lambda(S)$.

Since $\Y \cdot \Z \subset I_Q$, by Proposition \ref{MainProp},
$R'$ and $S'$ are Gorenstein Artin and $Q \simeq R' \#_\sk S'$.
Hence $\lambda(R') + \lambda(S') = \lambda(Q) +2 = \lambda(R) + \lambda(S)$. Since $\lambda(R') \leq \lambda(R)$
and $\lambda(S') \leq \lambda(S)$, we get $\lambda(R') = \lambda(R)$ and $\lambda(S') = \lambda(S)$, completing the proof.
\end{proof}

\subsection{A criterion in terms of socles}
In Remark \ref{CS}(e), we observed that if $Q$, $R$, and $S$ are Gorenstein Artin local rings with common residue field $\sk$ such that $Q \simeq R \#_\sk S$, then $\overline{Q} \simeq \overline{R} \times_{\sk}\overline{S}$, where $\ \bar{ }\ $ 
denotes going modulo the respective socles. A natural question is whether the converse is true. In the following proposition, we prove a partial converse in the graded case, by showing that $Q$ can be decomposed as a connected sum. However, we are unable to conclude anything about the components.

\begin{proposition}\label{Socle}
Let $Q$, $R$, and $S$ be graded Gorenstein Artin $\sk$-algebras with $\ll(Q) = \ll(R) = \ll(S) \geq 3$. Suppose $\overline{Q} \simeq \overline{R} \times_{\sk}\overline{S}$ as graded rings, where $\ \bar{ }\ $ denotes going modulo the respective socles. Then $Q$ can be decomposed as a connected sum over $\sk$. 
\end{proposition}
\begin{proof}

Let $P = R \times_\sk S$, $\pi: P \rightarrow \bar R \times_\sk \bar S$ be the natural projection, $\overline \varphi: \overline{R} \times_{\sk}\overline{S}   \rightarrow \overline{Q}$ be the given isomorphism, and set $\varphi = \overline \varphi \pi$. 

Write $\m_{\overline R}=\langle y_i\rangle_{1\leq i\leq m}$ and $\m_{\ov S}=\langle z_j\rangle_{1\leq j\leq n}$, where $\deg(y_i) = 1 = \deg(z_j)$. Choose $x_{1i}$, $x_{2j} \in \m_Q$ such that $\varphi(y_i)= \overline{x}_{1i}$ and $\varphi(z_j)= \overline{x}_{2j}$. Since $\ll(R),\ll(S)\geq 2$, we get $\edim(P)=\edim(\bar R) + \edim(\bar S) = \edim(\overline{Q})$. Hence $\m_{\overline Q}$ is minimally generated by $\{\overline{x}_{1i}, \overline{x}_{2j}:1\leq i \leq m, 1 \leq j \leq n\}$. Thus, $\soc(Q) \subset \m_Q^2$ forces $\m_Q = \langle x_{1i}, x_{2j}: 1\leq i \leq m, 1 \leq j \leq n\rangle$ with $\edim(Q) = m + n$.

Since $y_iz_j=0$, we get $x_{1i}x_{2j} \in \soc(Q)$ for each $i$, $j$. Since $\deg(x_{1i}) = 1 = \deg(x_{2j})$, and $\ll(Q) \geq 3$, we see that $\soc(Q) \neq \langle x_{1i}x_{2j}\rangle$ for any $i$, $j$, forcing $x_{1i}x_{2j}=0$ for each $i$ and $j$. Therefore, $Q$ can be decomposed as a connected sum over $\sk$ by Theorem \ref{MainThm}.    
\end{proof}

\subsection{Quotients of Decomposables}

If $(Q,\m_Q,\sk)$ is Gorenstein Artin and $0 \neq f \in \m_Q$, then $\overline Q = Q/(0:_Q f)$ is also Gorenstein. Suppose $Q$ can be decomposed as a connected sum. A natural question is whether $\overline Q$ can also be decomposed as a connected sum. The following is a condition which follows immediately from Theorem \ref{MainThm}.

\begin{proposition}\label{Quot1}
Let $(Q,\m_Q,\sk)$ be an Gorenstein Artin local ring which is decomposable as a connected sum of rings with embedding dimensions $m$ and $n$. Let $0 \neq f \in \m_Q$, and $L = \left( \frac{(0:_Q f) + \m_Q^2}{\m_Q^2}\right)$.
If $\lambda(L) < \min\{m,n\}$, then $Q/(0:_Q f)$ is also decomposable as a connected sum over $\sk$.

In particular, the conclusion holds if $(0:_Q f) \subset \m_Q^2$.
\end{proposition}

\begin{proof}
Let $\overline Q = Q/(0:_Q f)$. 
By the hypothesis on $Q$, we get $\m_Q = \langle y_1, \ldots, y_m, z_1, \ldots, z_n \rangle$, with $\y \cdot \z = 0$. Let $\ \bar{}\ $ denote going modulo $(0:_Q f)$. 

Since $\lambda(L) < \min\{m,n\}$, $\m_{\overline Q}$ is minimally generated by $\bar{y}_1, \ldots, \bar{y}_r, \bar{z}_1, \ldots, \bar{z}_s$ for some $r \geq 1$, and $s \geq 1$. Hence, by Theorem \ref{MainThm}, $\bar{\y} \cdot \bar{\z} = 0$ forces $\overline Q$ to be decomposable as a connected sum over $\sk$, of rings with embedding dimensions $r$ and $s$.
\end{proof}

In the next proposition, we give an answer to the question: When is $(0:_Q f) \subset \m_Q^2$?

\begin{proposition}\label{Quot2}
Let $Q$ and $f$ be as in Proposition \ref{Quot1}, $\m_Q = \langle \y, \z\rangle$ with $\y \cdot \z = 0$, $\soc(Q) = \langle \delta_Q \rangle$, and $f = f_y + f_z$, where $f_y \in \langle \y \rangle$, and $f_z \in \langle \z \rangle$. 

Then, $(0:_{Q} f)\subset \m^2_{Q}$ if and only if $(\langle\delta_Q\rangle :_Q f_y) \cap \langle \y \rangle \subset \langle \y \rangle^2$ and $(\langle\delta_Q\rangle :_Q f_z) \cap \langle \z \rangle\subset \langle \z \rangle^2$.  
\end{proposition}
\begin{proof}
($\Leftarrow$:)
Let $g \in (0:_Q f)$, and write $g = g_y + g_z$, for $g_y \in \langle \y \rangle$, and $g_z \in \langle \z \rangle$. Since $\y \cdot \z = 0 = fg$, we get $f_yg_y + f_zg_z = 0$. 
Hence $f_yg_y \in \langle \y \rangle \cap \langle \z \rangle = \soc(Q)$, where the last equality is by Remark \ref{CS}(h). Thus $g_y \in (\langle \delta_Q \rangle:_Q f_y) \cap \langle \y \rangle \subset \langle \y \rangle^2 \subset \m_Q^2$. By symmetry, $g_z \in \m_Q^2$, and hence $g \in \m_Q^2$.

($\Rightarrow$:) We prove $(\langle\delta_Q\rangle :_Q f_y) \cap \langle \y \rangle \subset \langle \y \rangle^2$. The proof of $(\langle\delta_Q\rangle :_Q f_z) \cap \langle \z \rangle\subset \langle \z \rangle^2$ is similar. 
Let $g_y \in (\langle\delta_Q\rangle :_Q f_y) \cap \langle \y \rangle$. If $f_yg_y = 0$, then $g_y \in (0:_Qf) \cap \langle \y \rangle \subset \m_Q^2 \cap \langle \y \rangle = \langle \y \rangle^2$. Hence, assume that $\langle f_yg_y \rangle = \soc(Q)$. Note that $\z \cdot f_y = 0$, and $(0:_{Q} f)\subset \m^2_{Q}$, force $f_z \neq 0$. Hence there is a $g_z \in \m_Q$ such that $f_z (- g_z) = f_yg_y$. Without loss of generality, we may assume that $g_z \in \langle \z \rangle$. 

Let $g = g_y + g_z$. Then $fg = f_yg_y + f_zg_z = 0$, hence $g \in 0:_Qf \subset \m_Q^2$. Hence, if $g_y \in \langle \y \rangle \setminus \m_Q^2$, then $g_z \in \langle \z \rangle \setminus \m_Q^2$, since $g \in \m_Q^2$. But then $\m_Q$ is minimally generated by a set containing $g_y$, and $g_z$, which contradicts $g_y + g_z \in \m_Q^2$. Thus, $g_y \in \m_Q^2$. Hence $g_y \in \langle \y \rangle \cap \m_Q^2 = \langle \y \rangle^2$. 
\end{proof}

\subsection{The equicharacteristic case}

Propositions \ref{CP2}, \ref{MainProp} and \ref{Prop2} give relations between the Cohen presentations of fibre products and connected sums with those of their components. In the equicharacteristic case, assuming all the rings are Artinian, we can write them as quotients of polynomial rings over the residue field. This allows us to give other relations among the defining ideals, which can be used in computations using computer algebra packages. As Remarks \ref{kalg} and \ref{Prop3} show, in the equicharacteristic case, the components of the decomposition, into either a fibre product or a connected sum, can be identified with subrings of the given ring.

We use the following notation in this subsection: 
For an ideal $J$ of $\sk[\Y]$ or $\sk[\Z]$, $J^e$ denotes its extension to $\sk[\Y,\Z]$ via the natural inclusions
$\sk[\Y] \hookrightarrow \sk[\Y,\Z]$ and $\sk[\Z] \hookrightarrow \sk[\Y,\Z]$ respectively.

The next remark is an analogue of Proposition \ref{CP2} in the equicharacteristic case.
\begin{remark}[Fibre Products in Equicharacteristic]\label{kalg}{\rm Let $R = \sk[\Y]/I_R$ and $S = \sk[\Z]/I_S$ be $\sk$-algebras with
$I_R \subset \langle \Y\rangle^2$ and $I_S \subset \langle\Z\rangle^2$. 
Let $P \simeq R \times_\sk S$. Since 
$R\simeq \sk[\Y,\Z]/(I_R^e + \langle \Z \rangle)$ and $S\simeq \sk[\Y,\Z]/(I_S^e + \langle \Y \rangle)$, a proof similar to Proposition \ref{CP2}(b) shows that $P  \simeq \sk[\Y,\Z]/I_P$, where $I_P = I_R^e + I_S^e + \langle \Y \cdot \Z\rangle$. 

Furthermore, we claim that $I_R = I_P \cap \sk[\Y]$, and $I_S = I_P \cap \sk[\Z]$. 
By symmetry, it is enough to show that $I_P \cap \sk[\Y] \subset I_R$.
Let $F \in I_P \cap \sk[\Y]$.
Since $I_P = I_R^e + I_S^e + \langle \Y \cdot \Z\rangle$, we can write $F = F_1 + F_2$, where $F_1 \in I_R$ and every term
of $F_2$ is a multiple of some $Z_j$. But $F_2 = F-F_1\in \sk[\Y]$, hence $F_2 = 0$. 
}\end{remark}

\begin{remark}[Connected Sums in Equicharacteristic]\label{Prop3}{\rm Let $(Q,\m_Q,\sk)$ be a Gorenstein Artin local $\sk$-algebra.
\begin{enumerate}[a)]
\item If $Q \simeq R \#_\sk S$, where $R \simeq \sk[\Y]/I_R$ and $S \simeq \sk[\Z]/I_S$ are Gorenstein Artin local rings, then we can write $Q \simeq \sk[\Y,\Z]/I_Q$, where 
$I_Q = I_R^e + I_S^e+ \langle \Y \cdot \Z\rangle+ \langle \Delta_R-  \Delta_S\rangle$, where $\Delta_R \in \sk[\Y]$ and $\Delta_S \in \sk[\Z]$ are such that $\langle \Delta_R \rangle R = \soc(R)$ and $\langle \Delta_S \rangle S = \soc(S)$. This follows from Definition \ref{CSdef} and Remark \ref{kalg}. 

Furthermore, by setting $I_{R'} = I_Q \cap \sk[\Y]$ and $I_{S'} = I_Q \cap \sk[\Z]$, and imitating the proof of Proposition \ref{Prop2}, one can see that $I_R = I_Q \cap \sk[\Y]$ and $I_S = I_Q \cap \sk[\Z]$. 

\item Conversely, suppose $Q$ can be decomposed as a connected sum over $\sk$. We can write $\m_Q = \langle y_1,\ldots, y_m,z_1,\ldots, z_n\rangle$, $m$, $n \geq 1$, such that $\y \cdot \z = 0$, and $\langle \y \rangle \cap \langle \z \rangle = \soc(Q)$. Set $Q = \sk[\y,\z] \simeq \sk[\Y,\Z]/I_Q$, where the $Y_i$'s and $Z_j$'s are mapped onto the $y_i$'s and $z_j$'s respectively. Then $Q \simeq R \#_\sk S$, where, by the last part in (a), we get $R \simeq \sk[\Y]/(I_Q \cap \langle \Y \rangle)$ and $S \simeq \sk[\Z]/(I_Q \cap \langle \Z \rangle)$. 

\item In particular, if $Q$ is standard graded, then $R$ and $S$ can be assumed to be so. Moreover, given $\Y \cdot \Z \subset I_Q$, one can use a computer algebra package (e.g., the elimination package in Macaulay2) to compute $I_R$ and $I_S$, and get a decomposition of $Q$ as a connected sum over $\sk$.
\end{enumerate}
}\end{remark}

An application of Proposition \ref{Quot2} is its analogue in the equicharacteristic case. Let the notation be as in the above remark.
Note that every $g \in \m_Q$, can be written as $g_y + g_z$, for some $g_y \in \m_R$, and $g_z \in \m_S$. In particular, with notation as in Proposition \ref{Quot2}, we can assume that $f_y \in \m_R$, and $f_z \in \m_S$. 

Furthermore, from the above description, it is clear that $\m_R = \m_Q \cap \sk[\y]$, $\m_R^2 = \langle \y \rangle^2Q \cap \sk[\y]$, and $\soc(R) = \soc(Q) \cap \sk[\y]$. Furthermore, if $\soc(R) = \langle \delta_R \rangle$, then, since $R$ is a subring of $Q$, and they are both Gorenstein Artin, we get $\soc(Q) = \langle \delta_R \rangle Q$.
The corresponding statements for $S$ are also true.

\begin{theorem}\label{Quot3}
With notation as above, let $(R,\m_{R},\sk)$, $(S,\m_{S},\sk)$ and $(Q,\m_{Q},\sk)$ be Gorenstein Artin $\sk$-algebras, such that $Q\simeq R\#_{\sk}S$. Let $f_y\in \m_R$, $f_z\in \m_S$, and $f = f_y + f_z \in \m_Q$. Then $(0:_{Q} f)\subset \m^2_{Q}$ if and only if  $(\langle \delta_R \rangle:_{R}f_y)\subset \m_{R}^2$ and $(\langle \delta_S \rangle:_{S}f_z)\subset \m_{S}^2$, where $\langle \delta_R \rangle = \soc(R)$, and $\langle \delta_S \rangle = \soc(S)$.
\end{theorem}

\begin{proof} 
Assume $(\langle \delta_R \rangle:_{R}f_y)\subset \m_{R}^2$, and let $g  \in (\langle\delta_Q\rangle :_Q f_y)$. We can write $g = g_y + g_z$, where $g_y \in \m_R$, and $g_z \in \langle \z \rangle \setminus \m_R$. Now, $f_y g_z \in \langle \y \rangle \cap \langle \z \rangle = \soc(Q)$. Hence $f_yg \in \soc(Q)$ implies that $f_yg_y \in \soc(Q) \cap \sk[\y] = \soc(R)$. Hence, $g_y \in (\langle \delta_R \rangle :_R f_y) \subset \m_R^2 \subset \langle \y \rangle^2$. If we further assume that $g \in \langle \y \rangle$, then $g_z = 0$.
Hence, $g = g_y$, proving $(\langle\delta_Q\rangle :_Q f_y) \cap \langle \y \rangle \subset \langle \y \rangle^2$.

On the other hand, let us assume $(\langle\delta_Q\rangle :_Q f_y) \cap \langle \y \rangle \subset \langle \y \rangle^2$. Suppose $g \in (\langle \delta_R \rangle:_{R}f_y)$. Clearly, $g \in \langle \y \rangle Q$, and $gf_y \in \langle \delta_R \rangle Q = \soc(Q)$. Thus $g \in (\langle\delta_Q\rangle :_Q f_y) \cap \langle \y \rangle$. Hence, the hypothesis implies that $g \in \langle \y \rangle^2 \cap \sk[\y] = \m_R^2$. 

Similar statements are true for $S$, and hence, this theorem follows from Proposition \ref{Quot2}.
\end{proof}

\section{Uniqueness of Decomposition}\label{Unique}

In the previous section, we have given equivalent conditions for a Gorenstein Artin local ring to be decomposable as a connected sum over $\sk$. A natural question to ask is whether such a decomposition is unique. 

Before we state the precise question, we introduce some notation, and state some properties of fibre products and connected sums, which follow from the results in Section \ref{FC} by induction. 

\begin{remark}[$p$-fold fibre products and connected sums]\label{Unique1}\hfill{}\\
{\rm Let $(R_i,\m_{R_i},\sk)$, $i = 1,\ldots, p$, be Gorenstein Artin local rings. 
\begin{enumerate}[a)]
\item By the universal property of pullbacks, one can see that taking fibre products over $\sk$, is both associative, and commutative, i.e., $R_1\times_{\sk}(R_2\times_{\sk}R_3)\simeq (R_1\times_{\sk}R_2)\times_{\sk}R_3$ and $R_1\times_{\sk}R_2\simeq R_2\times_{\sk}R_1$.

Hence, by Remark \ref{CS}(g), it follows that taking connected sums are associative and commutative, i.e., $R_1\#_{\sk}(R_2\#_{\sk}R_3)\simeq (R_1\#_{\sk}R_2)\#_{\sk}R_3$ and $R_1\#_{\sk}R_2\simeq R_2\#_{\sk}R_1$.

\item Set $P=R_1\times_{\sk}\cdots\times_{\sk}R_p$, and $Q=R_1\#_{\sk}\ldots\#_{\sk}R_p$. By induction, and Remark \ref{FP}(f), $P$ is local with maximal ideal $\m_P \simeq \prod_{i=1}^p \m_{R_i}$. 
Therefore, $\edim(P)=\sum_{i=1}^p\edim(R_i)$, $\lambda(P)=\sum_{i=1}^p\lambda(R_i)-(p-1)$, and if $\ell\ell(R_i)\geq 1$, then $\soc(P)=\bigoplus_{i=1}^p\soc(R_i)$.

\item Let $\pi_R: P \rightarrow Q$ be the natural projection. 
Let $\ov \m_{R_i}$ be the image of $\m_{R_i}$ in $Q$ under the map $\m_{R_i} \hookrightarrow \m_P \overset{\pi_R}\longrightarrow\!\!\!\!\rightarrow \m_Q$, and for $a_i \in \m_{R_i}$, denote the corresponding image in $Q$ by $\bar a_i$. 

The maximal ideal of $Q$ is $\m_Q = \sum_{i = 1}^p \overline{\m}_{R_i}$. Furthermore, $\lambda(Q)=\sum_{i=1}^p\lambda(R_i)-2(p-1)$, and $\ker(\pi_R) \subset \soc(P)$. Thus if $\ll(R_i)\geq 2$, then $\edim(Q)=\edim(P) = \sum_{i=1}^p\edim(R_i)$. 

\item For $1 \leq i \leq p$, let $R_i' = R_1\#_\sk \cdots \#_{\sk} R_{i-1}\#_{\sk} R_{i+1}\#_\sk \cdots \#_{\sk} R_p$. Observe that for each $i$, $Q \simeq R_i \#_{\sk} R_i'$, $\m_Q = \ov \m_{R_i} + \ov \m_{R_i'}$, 
and $\ov \m_{R_i} \cap \ov \m_{R_i'} = \soc(Q)$. Hence, by Remarks \ref{CS}(e), and \ref{FP}(f), we see that $\m_{R_i}/\soc(R_i) \simeq \m_Q/\ov \m_{R_i'} \simeq \ov \m_{R_i}/ \soc(Q)$ as $P$-modules. 

Thus, if $\ll(R_i) \geq 2$ for each $i$ (and hence $\ll(Q) \geq 2$), then $\{a_{i1}, \ldots, a_{ir_i}\}$ is a minimal generating set for $\m_{R_i}$ if and only if $\{\bar a_{i1}, \ldots, \bar a_{ir_i}\}$ is a minimal generating set of $\ov \m_{R_i}$. 
 \end{enumerate}
}\end{remark}

With this notation, we are now ready to state the precise question regarding the uniqueness of the decomposition of a Gorenstein Artin local ring $Q$ as a connected sum over $\sk$. Note that for a meaningful answer, we need to make two assumptions. The first is that each of the components appearing in a connected sum decomposition of $Q$ should themselves be indecomposable as a connected sum over $\sk$. Furthermore, to avoid the trivial decomposition, one can assume that their Loewy lengths are at least 2.
 
\begin{question}\label{uniqueness} Let $(R_i,\m_{R_i},\sk)$, $i = 1,\ldots, p$, and $(S_j,\m_{S_j},\sk)$, $j = 1,\ldots, q$ be Gorenstein Artin with $\ll(R_i), \ll(S_j) \geq 2$, such that they are indecomposable as connected sums over $\sk$.\\ If $\varphi: S_1\#_{\sk}\ldots\#_{\sk}S_q \overset{\simeq}\longrightarrow Q = R_1\#_{\sk}\ldots\#_{\sk}R_p$ is an isomorphism, is it necessary that\\ {\rm (i)} $p = q$, and {\rm (ii)} there is a permutation $\sigma$ of $\{1,\ldots, p\}$ such that $S_i\simeq R_{\sigma(i)}$ for each $i$?
\end{question}

Note that if the answer is affirmative for every $\varphi$ in Question \ref{uniqueness}, then the decomposition of $Q$ into indecomposables is unique. Assuming that $\ll(R_i), \ll(S_j) \geq 3$ for each $i$ and $j$, we prove this in the graded case (see Theorem \ref{3.5}), and show that in general, if $\varphi$ can be lifted to an isomorphism of the corresponding fibre products, then the above question has a positive answer (see Theorem \ref{3.4}). 

We will use the following setup in Remark \ref{Unique2}, and Lemmas \ref{3.3} and \ref{3.35}.
\begin{setup}\label{setup}\hfill{}{\rm
\begin{itemize}
\item For $1 \leq i \leq p$, let $(R_i,\m_{R_i},\sk)$ be Gorenstein Artin local rings, $P = R_1 \times_{\sk} \cdots \times_{\sk} R_p$, and $Q = R_1 \#_{\sk} \cdots \#_{\sk} R_p$.

\item Let $(S_1,\m_{S_1},\sk)$ and $(S_2,\m_{S_2},\sk)$ be Gorenstein Artin such that $\m_{S_1}$ and $\m_{S_2}$ are minimally generated by $\{ y_1,\ldots, y_m \}$, and $\{z_1, \ldots, z_n \}$ respectively. 

\item Let $\varphi: S_1 \#_\sk S_2 \overset{\simeq}\longrightarrow R_1 \#_{\sk}\cdots \#_{\sk} R_p$ be an isomorphism. For $1 \leq k \leq m$, $1\leq l \leq n$, $1 \leq i \leq p$, let $u_{ik}, v_{il} \in \m_{R_i}$ be such that $\varphi(\bar y_k) = \sum_{i = 1}^p\bar{u}_{ik}$, and $\varphi(\bar z_l) = \sum_{i = 1}^p\bar{v}_{il}$. 
\end{itemize}
}\end{setup}

We begin with the following remark.

\begin{remark}\label{Unique2}{\rm With the notation as in Setup \ref{setup}, let $\varphi(\ov \m_{S_1}) \subset \ov \m_{R_i}$ for some $i$. 
As in Remark \ref{Unique1}(d), we see that $\varphi$ induces an isomorphism $\ov \varphi: \prod_{j=1}^2 \m_{S_j}/\soc(S_j) \rightarrow \prod_{j = 1}^p \m_{R_j}/\soc(R_j)$, such that $\ov \varphi\left( \m_{S_1}/\soc(S_1)\right) \subset \m_{R_i}/\soc(R_i)$. Hence, if the Loewy lengths of all the rings involved are at least $2$, then, by Remark \ref{Unique1}(d), a minimal generating set of $\m_{S_1}$ corresponds to a part of a minimal generating set of $\m_Q$ in $\ov \m_{R_i}$, and hence corresponds to a part of a minimal generating set of $\m_{R_i}$. In particular, $\edim(S_1) \leq \edim(R_i)$.
}\end{remark}

The following is a key lemma which follows from Proposition \ref{Prop2}.

\begin{lemma}\label{3.3}
Let the notation be as in Setup \ref{setup}, with $p = 2$, and $\ell\ell(R_i),\ell\ell(S_j)\geq 2$. If $\varphi({\overline{\m}_{S_1}})\subseteq \overline{\m}_{R_1}$ and $\varphi(\overline{\m}_{S_2})\subseteq \overline{\m}_{R_2}$, then $S_1\simeq R_1$ and $S_2\simeq R_2$.    
\end{lemma}
\begin{proof}

Set $Q=S_1\#_{\sk}S_2$, and let $\m_{S_1} = \langle \y \rangle$, and $\m_{S_2} = \langle \z \rangle$, with corresponding images $\overline \y$, and $\overline \z$ in $\m_Q$. Let $Q \simeq \widetilde Q/I_Q$ be a Cohen presentation of $Q$, with $\m_{\widetilde Q} = \langle \Y, \Z \rangle$, where the $Y$'s and $Z$'s are the respective preimages of the $\overline y$'s and $\overline z$'s respectively. By Proposition \ref{Prop2} $S_1 \simeq \widetilde Q/J_{S_1}$, where $J_{S_1} = (I_Q \cap \langle \Y \rangle) + \langle \Z \rangle$. We show that $S_1 \simeq R_1$, the other isomorphism follows similarly.

Since $\varphi(\overline{\m}_{S_i})\subseteq\overline{\m}_{R_i}$, by Remark \ref{Unique2}(c), we have $\edim(S_i)\leq \edim(R_i)$, for $i = 1,2$. The equalities are forced since $\edim(R_1)+\edim(R_2)=\edim(S_1)+\edim(S_2)$. In particular, $\overline \m_{R_1}$ and $\overline \m_{R_2}$ are minimally generated by $\varphi(\y)$, and $\varphi(\z)$ respectively, and $\varphi(\y)\cdot \varphi(\z) = 0$. Hence, by Proposition \ref{Prop2}, $R_1 \simeq \widetilde Q/J_{R_1}$, where $J_{R_1} = (I_Q \cap \langle \Y \rangle) + \langle \Z \rangle$. This proves the lemma.
\end{proof}

A consequence is the following technical lemma. Certain results related to the uniqueness of decomposition follow from this as an application, e.g., see Theorems \ref{3.4} and \ref{3.5}.

\begin{lemma}\label{3.35}
Let the notation be as in Setup \ref{setup}. If each $R_i$ is indecomposable as a connected sum over $\sk$, and $u_{ik}v_{il} = 0$ for each $i$, $k$, and $l$, then there is a permutation $\sigma$ of $\{1,\ldots, p\}$, and $1 \leq t < p$, such that $S_1\simeq R_{\sigma(1)}\#_{\sk} \cdots \#_{\sk} R_{\sigma(t)}$ and $S_2\simeq R_{\sigma(t+1)}\#_{\sk} \cdots \#_{\sk} R_{\sigma(p)}$.

In particular, if $S_1$ is indecomposable as a connected sum over $\sk$, then $t = 1$,\\ i.e., $S_1 \simeq R_{\sigma(1)}$, and $S_2 \simeq R_{\sigma(2)}\#_{\sk} \cdots \#_{\sk} R_{\sigma(p)}$.
\end{lemma}

\begin{proof}
Let $Q = R_1 \#_{\sk}\cdots \#_{\sk} R_p$. Firstly note that since $\varphi$ is an isomorphism, $\{\varphi(\bar y_k), \varphi(\bar z_l)\}$ is a minimal generating set for $\m_Q$. Hence, by Remark \ref{Unique1}(d), one can see that the set $\{u_{ik}, v_{il}: 1 \leq k \leq m, 1 \leq l \leq n\}$ generates $\m_{R_i}$, for each $i$. Hence, by Theorem \ref{MainThm}, the indecomposibility of $R_i$ as  a connected sum over $\sk$, and the fact $u_{ik}v_{il} = 0$ for each $k$ and $l$, forces $\m_{R_i} =  \langle u_{i1},\ldots,u_{im}\rangle$ or $\m_{R_i} =  \langle v_{i1},\ldots,v_{in}\rangle$.

Since $\{\varphi(\bar y_k), \varphi(\bar z_l)\}$ is a minimal generating set for $\m_Q$, there exist $i, j \geq 1$, $i \neq j$ such that $\m_{R_i} =  \langle u_{i1},\ldots,u_{im}\rangle$, and $\m_{R_j} =  \langle v_{j1},\ldots,v_{jn}\rangle$. Thus, there exists $t \in \{1, \ldots, p-1\}$, $A = \{i_1,\ldots, i_t\}$, $B = \{j_{1}, \ldots, j_{p-t}\}$ with $A \cap B = \phi$, and $A \cup B = \{1,\ldots,p\}$, such that $\m_{R_i} =  \langle u_{i1},\ldots,u_{im}\rangle$ for $i \in A$, and $\m_{R_j} =  \langle v_{j1},\ldots,v_{jn}\rangle$ for $j \in B$.

Let $i \in A$, and $l \in \{1,\ldots, n\}$. By Theorem \ref{MainThm}, $\bar v_{il}\ov \m_{R_j} = 0$ for $j \neq i$. Furthermore, since $\m_{R_i} =  \langle u_{i1},\ldots,u_{im}\rangle$, $u_{ik}v_{il} = 0$ implies that $v_{il} \ov \m_{R_i} = 0$. Thus, $\m_Q = \sum_{j = 1}^p \ov \m_{R_j}$ forces $\bar v_{il} \in \soc(Q) \subset \ov \m_{R_j}$ for any $j \in B$. Thus, $\varphi(\bar z_l) \in \ov \m_{R_{j_1} \#_\sk \cdots \#_\sk R_{j_{p-t}}}$ for each $l$. 

Similarly, we can see that $\varphi(\bar y_k) \in \ov \m_{R_{i_1} \#_\sk \cdots \#_\sk R_{i_t}}$ for each $k$.
Thus, by Lemma \ref{3.3}, we get $S_1\simeq R_{i_1}\#_{\sk}\ldots \#_{\sk}R_{i_t}$ and  $S_2\simeq R_{j_1}\#_{\sk}\ldots \#_{\sk}R_{j_{p-t}}$, and if $S_1$ is indecomposable, then $t = 1$.

The proof is complete by letting $\sigma$ be the permutation of $\{1, \ldots, p\}$ given by $\sigma(k) = i_k$ for $1 \leq k \leq t$, and $\sigma(k) = j_{k-t}$ for $t+1 \leq k \leq p$.
\end{proof}

In the following theorem, we give two scenarios where Question \ref{uniqueness} has a positive answer.

\begin{theorem}\label{3.4}
Let $Q$, $R_i$, $S_j$ and $\varphi$ be as in Question \ref{uniqueness}. Then 
{\rm{(i)}} $p=q$, and {\rm{(ii)}}  there is a permutation $\sigma$ of $\{1,\ldots, p\}$ such that $S_i\simeq R_{\sigma(i)}$ for each $i$ if either one of the following conditions are satisfied.\\
{\rm(a)} $\varphi$ lifts 
to an isomorphism $\tilde \varphi : S_1\times_{\sk}\cdots\times_{\sk}S_q \overset{\simeq}\longrightarrow R_1\times_{\sk}\cdots\times_{\sk}R_p$, i.e., the following diagram commutes: \[
 	\xymatrixrowsep{8mm} \xymatrixcolsep{8mm}
 	\xymatrix{
		S_1\times_{\sk}\ldots\times_{\sk}S_q \ar[d]_{\pi_{S}}\ar[r]^{\tilde \varphi}& R_1\times_{\sk}\ldots\times_{\sk}R_p \ar[d]^{\pi_{R}}\\
 		S_1\#_{\sk}\ldots\#_{\sk}S_q \ar[r]^{\varphi}& R_1\#_{\sk}\ldots\#_{\sk}R_p
	}
	\]where $\pi_R$ and $\pi_S$ are the natural surjective maps. \\  
{\rm(b)} $R_i$, and $S_j$ are graded $\sk$-algebras for each $i$ and $j$, and $\varphi$ is a graded homomorphism.
\end{theorem}
\begin{proof}
Without loss of generality, we may assume that $2 \leq q \leq p$. Set ${Q_i}=S_i\#_{\sk}\ldots\#_{\sk}S_q$. Note that $S_{i}\#_{\sk}{Q_{i+1}} \simeq {Q_{i}}$ and $S_1\#_{\sk} {Q_2} \overset{\varphi}\simeq Q$. 
Let $\m_{S_1}$ and $\m_{Q_2}$ be minimally generated by $\{ y_1,\ldots, y_m\}$, and $\{z_1, \ldots, z_n\}$, respectively. For $1 \leq k \leq m$, $1\leq l \leq n$, $1 \leq i \leq p$, let $u_{ik}, v_{il} \in \m_{R_i}$ be such that $\varphi(\bar y_k) = \sum_{i = 1}^p\bar{u}_{ik}$, and $\varphi(\bar z_l) = \sum_{i = 1}^p\bar{v}_{il}$. We first claim $u_{ik}v_{il} = 0$ for each $i$, $k$, and $l$, assuming either hypothesis (a) or (b).

Suppose hypothesis (a) is true. Let $z_l' \in S_2 \times_\sk \cdots \times_{\sk} S_q$ be the corresponding preimages of the $z_l$'s. Since $\varphi \pi_S = \pi_R \tilde \varphi$, for each $i$, we can assume that $u_{ik}, v_{il} \in \m_{R_i}$ are such that $\tilde \varphi(y_k) = \sum_{i = 1}^p{u}_{ik}$ and $\tilde \varphi( z_l') = \sum_{i = 1}^p{v}_{il}$ for each $k$ and $l$. 

Let $P = R_1 \times_{\sk} \cdots \times_{\sk} R_p$. Now, $y_k  z_l' = 0$ forces $\sum_{i = 1}^p{u}_{ik} v_{il} = 0$ in $\m_P$. Since $\m_P \simeq \prod_{i=1}^p \m_{R_i}$, we see that $u_{ik}v_{il} = 0$ for each $i$, $k$, and $l$. 

Assume hypothesis (b). Since $\varphi$ is graded, we see that for each $i$ and $k$, either $\deg(u_{ik}) = 1$, or $u_{ik} = 0$. The same holds for $v_{il}$ for each $i$ and $l$. Fix $k$ and $l$. Since $\bar y_k \bar z_l = 0$, we get $\sum_{i = 1}^p \bar u_{ik}\bar v_{il} = 0$. In particular, $\sum_{i = 1}^p u_{ik}v_{il} \in \ker(\pi_R) \subset \soc(P) = \prod_{i=1}^p \soc(R_i)$. Thus, we get $u_{ik}v_{il} \in \soc(R_i)$ for each $i$. But $\ll(R_i) \geq 3$, hence degree considerations force $u_{ik}v_{il} = 0$, for each $i$, $k$, and $l$.

Thus, assuming either (a) or (b), we see that $u_{ik}v_{il} = 0$ for each $i$, $k$, and $l$, where $\varphi(\bar y_k) = \sum_{i = 1}^p\bar{u}_{ik}$, and $\varphi(\bar z_l) = \sum_{i = 1}^p\bar{v}_{il}$. Hence, by Lemma \ref{3.35}, we get $S_1 \simeq R_{i_1}$ for some $1 \leq i_1 \leq p$, and $Q_2 \simeq R_{i_2} \#_\sk \cdots \#_{\sk} R_{i_p}$, where $\{i_2,, \ldots, i_p\} = \{1,\ldots,p\} \setminus \{i_1\}$. 

We finish the proof of the theorem by induction on $q$. If $q = 2$, then $S_1 \simeq R_{i_1}$, and $S_2 \simeq Q_2$ forces $Q_2$ to be indecomposable. Since $\ll(R_i) \geq 2$, this is only possible if $p = 2$. Thus $S_2 \simeq R_{i_2}$, where $\{i_1, i_2\} = \{1,2\}$, proving the result in this case.

If $q \geq 3$, then, by induction, $S_2\#_{\sk} \cdots \#_{\sk} S_q = Q_2 \simeq R_{i_2} \#_\sk \cdots \#_{\sk} R_{i_p}$ gives $p-1 = q-1$, i.e., $p = q$. Furthermore, for each $i \in \{2, \ldots, q =p\}$, there is a rearrangement $\{\sigma(2), \ldots, \sigma(p)\}$ of $\{i_2,\ldots, i_p\}$ such that $S_i \simeq R_{\sigma(i)}$. Letting $\sigma(1) = i_1$, the proof is complete.
\end{proof}

Thus, Theorem \ref{3.4} answers Question \ref{uniqueness} in the affirmative for every isomorphism in the graded case. This proves the uniqueness of decomposition of a graded Gorenstein Artin $\sk$-algebra into graded indecomposables, which we record below.

\begin{theorem}[Uniqueness of Decomposition as Connected Sums: The Graded Case]\label{3.5}\hfill{}\\
Let $Q$ be a graded Gorenstein Artin $\sk$-algebra with $\ll(Q) \geq 3$. Then, the decomposition of $Q$ into graded components, which are indecomposables as a connected sum over $\sk$, is unique up to isomorphism. 
\end{theorem} 
\begin{proof}
Let $R_i$, and $S_j$ {\rm(}$1 \leq i \leq p$, and $1 \leq j \leq q${\rm)}, be graded Gorenstein Artin $\sk$-algebras with $\ll(R_i) = \ll(S_j) \geq 3$, where each $R_i$ and $S_j$ is indecomposable as a connected sum over $\sk$. If $Q \simeq R_1\#_{\sk}\ldots\#_{\sk}R_p \simeq S_1\#_{\sk}\ldots\#_{\sk}S_q$ as graded rings, then by Theorem \ref{3.4}(b), we get $p=q$. Furthermore, there is a permutation $\sigma$ of $\{1,\ldots, p\}$ such that $R_{\sigma(i)}\simeq S_i$ for $1\leq i\leq p$,  i.e., the two decompositions of $Q$ are unique up to isomorphism.
\end{proof}

We end this article with two questions related to the uniqueness of decomposition in the general case. If either of the questions has a positive answer, then the decomposition of any Gorenstein Artin local ring into indecomposables is unique.

\begin{question}
Let $R_i$, $S_j$, $\varphi$ be as in Question \ref{uniqueness}.\\ {\rm i)} Does $\varphi$ lift to an isomorphism $\tilde \varphi: S_1 \times_\sk \cdots \times_{\sk} S_q \rightarrow R_1 \times_{\sk} \cdots \times_{\sk} R_p$?\\
{\rm ii)} Can $\varphi$ be modified to construct an isomorphism $\psi: S_1 \#_\sk \cdots \#_{\sk} S_q \rightarrow R_1 \#_{\sk} \cdots \#_{\sk} R_p$, which can be lifted to an isomorhism $\tilde \psi: S_1 \times_\sk \cdots \times_{\sk} S_q \rightarrow R_1 \times_{\sk} \cdots \times_{\sk} R_p$?
\end{question}

\noindent
{\bf Acknowledgement}\\
We would like to thank L. L. Avramov for comments inspiring this work. We would also like to thank him, 
A. A. Iarrobino and J. Weyman for insightful discussions.

\end{document}